\documentclass[11pt]{amsart}
\usepackage{amssymb, verbatim}

\newcommand{\FT}{{\mathcal{F}}}
\newcommand{\h}{{\frac 12}}
\newcommand{\pa}{{\partial}}

\newcommand{\RR}{{\mathbb{R}}}

\newcommand{\WF}{{\operatorname{WF}}}

\newcommand{\CI}{{\mathcal{C}^\infty}}
\newcommand{\CIc}{{\mathcal{C}_c^\infty}}
\newcommand{\abs}[1]{{\left\lvert{#1}\right\rvert}}

\newcommand{\Lap}{{\Delta}}

\newcommand{\schwartz}{{\mathcal{S}}}

\newcommand{\QSC}{{\mathrm{qsc}}}
\newcommand{\WFqsc}{{\WF_{\QSC}}}
\newcommand{\SC}{{\mathrm{sc}}}
\newcommand{\WFsc}{{\WF_{\SC}}}

\DeclareMathOperator{\Ai}{{Ai}}

\DeclareMathOperator{\supp}{{supp}}

\theoremstyle{plain}
\newtheorem{theorem}{Theorem}
\newtheorem{proposition}[theorem]{Proposition}

\theoremstyle{definition}
\newtheorem{definition}[theorem]{Definition}

\theoremstyle{remark}
\newtheorem*{remark}{Remark}

\author{Andrew Hassell}
\author{Jared Wunsch}

\title{On the structure of the Schr\"odinger propagator}

\begin{document}
\maketitle
\begin{abstract}
We discuss the form of the propagator $U(t)$ for the time-dependent
Schr\"odinger equation on an asyptotically Euclidean, or, more generally,
asymptotically conic, manifold with no trapped geodesics.  In the
asymptotically Euclidean case, if $\chi \in \mathcal{C}_0^\infty$, and with
$\mathcal{F}$ denoting Fourier transform, $\mathcal{F}\circ e^{-ir^2/2t}
U(t) \chi$ is a Fourier integral operator for $t\neq 0.$ The canonical
relation of this operator is a ``sojourn relation'' associated to the
long-time geodesic flow.  This description of the propagator follows from
its more precise characterization as a ``scattering fibered Legendrian,''
given by the authors in a previous paper and sketched here.

A corollary is a propagation of singularities theorem that permits a
complete description of the wavefront set of a solution to the
Schr\"odinger equation, restricted to any fixed nonzero time, in terms of
the oscillatory behavior of its initial data.  We discuss two examples
which illustrate some extremes of this propagation behavior.
\end{abstract}

\section{Introduction}
Let us consider the Schr\"odinger initial value problem for a particle moving in curved
space with metric $g$
\begin{equation}
\left(D_t + \h \Lap+V\right)\psi(z,t)=0, \quad \psi(0, \cdot) = \psi_0,
\label{scheqn}
\end{equation}
with $D_t=-i\pa_t$ and $\Lap$ the nonnegative Laplace-Beltrami operator
$$
-\frac 1{\sqrt g} \pa_{z^i} g^{ij}\sqrt g \pa_{z^j},
$$
and where $V$ is a smooth, real-valued potential function. The solution is given in terms of functional calculus by the formula
$$
\psi(t, \cdot) = e^{-it H} \psi_0, \quad H = \h \Lap + V,
$$
and the operator $e^{-itH}$ (or its kernel) is called the \emph{propagator}, or fundamental solution. 
On flat $\RR^n$ with $V=0,$ there is an explicit formula for the propagator
for \eqref{scheqn}: its Schwartz kernel is
\begin{equation}
K(z,w,t)=(2\pi i t)^{-n/2} e^{i\abs{z-w}^2/2t}.
\label{euclidian.propagator}
\end{equation} 
This may also be thought of as the solution of \eqref{scheqn} with initial
data the delta function at $w$ (the initial value problem makes sense for
any tempered distribution as initial data).

This solution exhibits some peculiar properties, from the point of view of
propagation phenomena.  Fixing $w=w_0 \in \RR^n$, we find that for any
$t>0,$ $K(t,z,w_0)$ is \emph{smooth}. Hence the initial delta-function
singularity at $z=w_0$ apparently disappears.  Additionally, since the
propagator is, by self-adjointness of $\Lap,$ unitary on $L^2,$ we can
reverse this process: If we take
$$
\psi_0=e^{i\lambda \abs{z}^2/2+\xi \cdot z}
$$
as smooth initial data for a solution to \eqref{scheqn}, then at
$t=-\lambda^{-1}$ the solution develops a delta-function singularity,
located at $z=-\xi/\lambda.$  One can see this either by explicitly
convolving with $K(z,w,t)$ and doing a Gaussian integral, or by noting
that this initial data is nothing but a multiple of a time-translated
version of the fundamental solution.

The upshot, then, is that singularities to solutions of the Schr\"odinger
equation can both appear and disappear.  In this paper, we address these
phenomena not just on $\RR^n,$ but more generally on curved space.  The
basic questions are where do singularities go when they disappear, and,
conversely, what causes their appearance?  (These questions are based on
the idea that, in evolution described by a unitary group, things should not
``disappear'' and ``appear,'' but rather be transformed into something more or
less equivalent.)  The most powerful tool for answering propagation
questions is to understand the kernel of the propagator: a sufficiently
precise understanding of it, particularly its asymptotics at infinity, will
allow us to answer both questions.  We describe a generalization of the
formula \eqref{euclidian.propagator} and an interpretation thereof which
enables us to describe the formation and disappearance of singularities.
We should emphasize that by a singularity at time $t$ we mean a point in
the \emph{wavefront set} of $\psi(t, \cdot)$ in the sense of H\"ormander
\cite{Hormander1}.  We recall that a point in the wavefront set of a
distribution is an element of the unit cosphere bundle, $(z, \hat\zeta)$
with $|\hat\zeta| = 1$, which intuitively describes an infinitesimal wave located at $z$ with wavefronts normal to $\hat\zeta$.

This note is intended as a ``user's guide'' to the authors' more technical
paper \cite{Hassell-Wunsch1} which describes both a parametrix construction
for the Schr\"odinger propagator and, as a corollary, a full propagation
of singularities theorem answering the two questions posed above. 

\

Much has long been known about the structure of the propagator on $\RR^n$
with nonzero potentials $V;$ indeed a number of rather precise parametrix
constructions exist (among others, those of Fujiwara \cite{Fujiwara1},
Zelditch \cite{Zelditch2}, Tr\`eves \cite{Treves1}, Yajima \cite{Yajima1}).

On curved space, on the other hand, almost nothing was known about the
propagator until comparatively recently.  The first variable-coefficient
results were those of Kapitanski-Safarov \cite{KapSaf1}, who showed that in
the case of a compactly-supported metric perturbation of $\RR^n$ with
\emph{no trapped geodesics} the kernel of the propagator is smooth for all
$t>0.$ The same authors subsequently constructed a parametrix
\cite{KapSaf2}, albeit without the control at infinity which will turn out
to be essential for the purposes at hand.  At around the same time,
Craig-Kappeler-Strauss proved the first result about \emph{microlocal}
regularity of solutions to \eqref{scheqn}.  Loosely speaking, they showed
that on asymptotically Euclidean space, if the initial data $\psi_0$ is a
Schwartz function in a conic microlocal neighborhood of a backward geodesic
passing through a point $(z,\hat \zeta) \in S^* \RR^n,$ then $(z,\hat
\zeta) \notin \WF \psi(t)$ for all $t>0.$ In particular, let $\chi(z)$
denote a cutoff function supported in a positive cone near infinity.
Assume that $\chi(z) \psi_0 \in \schwartz.$ Now let $\gamma(s)$ be any
geodesic so that $\gamma(s) \in \supp \chi$ for all $s<<0.$ We conclude
that $(\gamma(s), g_{ij} (\gamma')^j(s)) \notin \WF \psi(t,z)$ for any
fixed $t>0.$ Thus a single hypothesis gives regularity along the
(co-)tangents to a whole \emph{pencil} of geodesics emanating from $\supp
\chi$ for all $t>0.$ In Euclidean space, this yields regularity along the
tangents to all lines lying in a pencil of directions, and thus constrains
the direction in which singularities of $\psi$ can lie.

The second author \cite{Wunsch1} refined the result of
Craig-Kappeler-Strauss by introducing the \emph{quadratic wavefront set} ``at
infinity,'' measuring the quadratic oscillation of a distribution\footnote{
For example, in the case of a distribution $e^{i \phi(z)},$ where $\phi$ is
homogeneous of degree two, hence determined by a function $\tilde \phi$ on $S^{n-1}$,  the quadratic wavefront set is essentially the graph of $\tilde \phi$ over the sphere at infinity.}. It is analogous to the \emph{scattering wavefront set} previously introduced by Melrose
\cite{MR95k:58168} to describe linear oscillations.  In \cite{Wunsch1}
conditions were given, in terms of the quadratic scattering wavefront set, that constrain not only the
directions in which singularities can appear, but also the \emph{times} at
which they can appear.  The actual locations, however, remain mysterious
from this point of view. Theorems such as these, which determine $\zeta$
and $t$ but not $z$ in the wavefront set $(z, \zeta) \in \WF \psi(t)$ have
been called \emph{microglobal results} \cite{Steve}.

The quadratic-scattering wavefront set results of \cite{Wunsch1} have
recently been extended, using a rather different set of tools, to the
analytic category by Robbiano-Zuily \cite{Robbiano-Zuily}; in this
setting, even defining the quadratic scattering wavefront set involves an
appropriate two-parameter version of the FBI transform, rather than
the pseudodifferential methods employed in
\cite{MR95k:58168,Wunsch1}.

The first author thanks the Australian Research Council for its support, through a Fellowship and a Linkage grant, of this research. 
The second author acknowledges the support of the National Science
Foundation, through grants DMS-0100501 and DMS-0323021.

\section{The geometry}
The results of \cite{Hassell-Wunsch1} hold for a rather general class of
manifolds with large conic ends, known as \emph{scattering manifolds},
introduced by Melrose \cite{MR95k:58168}.  However for simplicity we will
restrict ourselves here to a sub-class consisting of asymptotically
Euclidean spaces.  Let $z$ be a Euclidean coordinate and let
$$
\theta=\frac z{\abs z}, r=\abs z
$$
be polar coordinates.  We will assume that our metric is smooth and has the form
\begin{equation}
g=\sum_{i=1}^n (dz^i)^2+ h,\quad h =  \frac{m \, dr^2}{r} + \frac{k_{ij}(r^{-1},\theta)dz^i dz^j}{r^2} \text{ for } r \gg 0,
\label{metric}\end{equation}
where $m$ is a constant and $k_{ij}$ is a $\CI$ function of its arguments, i.e.\ has an
asymptotic expansion (Taylor series) in descending powers of $r$.

We further make the following crucial nontrapping assumption\footnote{It is actually not necessary to make this assumption, but if we do not then the results only apply in the non-trapping part of phase space.} for $g$:
\begin{equation}
\text{For any geodesic }\gamma(t),\ \lim_{t\to\pm\infty} r(\gamma(t))=\infty.
\label{notrapping}
\end{equation}

As for the potential $V$, we assume $V \in \CI(\RR^n; \RR)$ and for $r\gg 0,$
\begin{equation}
V = \frac{c}{r} + \frac{\tilde V(r^{-1},\theta)}{r^2}
\label{pot-form}\end{equation}
with $c$ a constant and $\tilde V$ a $\CI$ function of its arguments.
If $m$ in \eqref{metric} and $c$ in \eqref{pot-form} both vanish, then $H$ is said to be \emph{short-range}, otherwise \emph{(gravitational) long-range}. 

We henceforth denote the Hamiltonian
$$
H=\h \Lap + V,
$$
hence the propagator is $e^{-it H}.$

\section{The form of the propagator}

In order to motivate our results about the Schwartz kernel of the
propagator, let us re-examine the special form
\eqref{euclidian.propagator} of the Euclidean one.  If we introduce polar
coordinates in the $w$ variable, we have
$$
K(t,r,\theta, z)= e^{i r^2/2t} \left( ae^{-i z\cdot \theta r/t}\right )
$$ where $a=(2\pi i t)^{-n/2} e^{iw^2/2t}.$ We have factored $K$ into this
form for several reasons.  First, we are considering $K$ as a function of
$r,\theta$ alone, with $t > 0$ and $z$ fixed.  Hence $a$ is, from this point of view,
constant.  Second, while $K$ is a rather uninteresting distribution from
the point of view of wavefront set, it is certainly not a Schwartz
function, owing to its oscillatory behavior at infinity.  We have chosen to
exhibit this oscillatory behavior by separating out the leading order
factor $e^{i r^2/2t},$ from the milder, plane wave oscillation of
$e^{-iz\cdot \theta r/t}.$ Recall now that part of our goal in constructing
the propagator is to understand the fate of the delta-function
singularity launched from the point $z \in \RR^n.$ From this point of view,
the leading order term $e^{ir^2/2t}$ is useless: it retains no
information about where the initial delta function lay, and merely
records, in its frequency $t^{-1}$, the elapsed time.  By contrast, the plane
wave oscillation term is of great interest: we can recover $z$ from
$e^{-iz\cdot \theta r/t}.$ Finally, note that $t$ appears in the phase in a
very simple way.

Since the quadratic term is of little interest, we will separate it explicitly in our formula for the propagator.  Thus, for
$t \neq 0,$ and for $H$ short-range, let
$$
W_t=e^{-ir^2/2t} e^{-itH}.
$$ In the Euclidean case, this is a constant multiple of a plane wave
$e^{-iz \cdot \theta r/t}.$ Fourier transforming in $w=r\theta,$ then,
gives a delta distribution $\delta(z-w)$.  This turns out to be a rather general phenomenon.

Let $\mathcal{F}$ denote the Fourier transform and let $\chi \in \CIc(\RR^n).$ One of the main results from \cite{Hassell-Wunsch1} can be expressed as
\begin{theorem}\label{par} Suppose $H$ is short-range and $t \neq 0$. Then the operator
$$
\mathcal{F}\circ W_t \chi
$$
is a classical Fourier integral operator of order $0$ on $\RR^n$ associated to a canonical transformation of $S^* \RR^n$. 
\label{thm:FIO}
\end{theorem}
\emph{Remark.} The point of multiplying on the right by $\chi$ is to localize the right
spatial variable to a compact set. Theorem~\ref{par} is essentially
concerned with the behaviour of the propagator when the left variable
approaches infinity. What happens as both variables independently approach
infinity is more complicated.

\vskip 5pt

Let $Q$ denote the Fourier integral operator of Theorem~\ref{par}, and $\gamma$ the contact transformation. Then a standard property of FIOs \cite{Hormander1} is that $Q$ moves wavefront set according to $\gamma$: 
$$
\WF (Q f) \subset \gamma(\WF f).
$$
It is now easier to see why the propagator creates and destroys wavefront
set: can can write $e^{-it H}\chi = e^{ir^2/2t} \FT^{-1} Q$.  Thus $Q$ moves wavefront set around, but $\FT^{-1}$
kills it (if it is compactly supported).  So a better way to think of the
above theorem is that $W_t$ maps wavefront set to ``Fourier transformed
wavefront set.''

In the gravitational long-range case, we need to modify the radial variable slightly.
Let
\begin{equation}
\tilde r = r + \frac{m}{2} \log r,
\label{rt}\end{equation}
and define
$$
W_t = e^{-i\tilde r^2/2t} e^{-itH}.
$$
This modification is familiar in, for example, the structure of generalized eigenfunctions or of Dollard wave operators for the Schr\"odinger operator with Coulomb potential \cite{Dollard}, \cite{Schiff}. 
Then $\mathcal{F}\circ W_t \chi$ is a quasi-classical\footnote{\emph{Quasi-classical} here means that the symbol $a(x, y, \theta)$ has an expansion that includes log terms as well as powers of $\theta$.} Fourier integral operator of order zero. 

To understand the canonical relation of Theorem~\ref{thm:FIO} will require
understanding the long-time limit of geodesic flow.  We begin, however, by
clarifying the action of $W_t$ itself on singularities of distributions.

\section{Scattering wavefront set}
We have seen that in order to describe the mapping properties of the
propagator on wavefront set, it is helpful to keep track of Fourier
transformed wavefront set.  The wavefront set of a tempered distribution on
$\RR^n$ is a closed subset of $S^*\RR^n;$ it is convenient for our purposes
to identify this space with the product $\RR^n \times S^{n-1}.$

Fourier transforming a compactly supported distribution with nonempty
wavefront set---a delta-function for instance---yields a smooth function
with oscillation at infinity---e.g.\ a plane wave.  In tracking where the
singularity has gone under Fourier transform, it is therefore helpful to
introduce a wavefront set that measures oscillation at infinity.  On
$\RR^n,$ we simply use the Fourier transform to do this.
\begin{definition}\label{defn:scwf}
The \emph{scattering wavefront set} of a distribution $u \in
\schwartz'(\RR^n)$ is the closed subset
$$
\WFsc u\subset S^{n-1}_{\hat z} \times \RR^n_\zeta
$$
given by
$$
(\hat z, \zeta) \in \WFsc u \Longleftrightarrow (\zeta,\hat z) \in
\WF \mathcal{F} (u).
$$
\end{definition}
So, for instance,
$$
\WFsc e^{i \xi \cdot z} =\{(\hat z, \xi): \hat z \in S^{n-1}\}
$$
captures the frequency of plane wave oscillation.

\begin{remark}
The scattering wavefront set was introduced by Melrose \cite{MR95k:58168}
in the more global context of scattering manifolds, where the definition is
more subtle.  Melrose's definition encompasses the set defined above, which
we think of as a subset of an appropriately scaled cotangent bundle at
infinity, and the ordinary wavefront set inside the cosphere bundle, as
well as a third component which interpolates between the two, and is a
subset, in our non-invariant notation, in $S^{n-1}_{\hat z} \times
S^{n-1}_{\hat \zeta}$ (the cosphere bundle at infinity).
\end{remark}

The \emph{quadratic-scattering wavefront set} of \cite{Wunsch1} can be
defined quite similarly: again a subset of $S^{n-1} \times\RR^n,$ it can be
defined by
$$
\WFqsc u = \WFsc \tilde u
$$
where $\tilde u$ is the distribution defined by
$$
\tilde u(z)=u(z/\sqrt{\abs{z}}).
$$

Theorem~\ref{thm:FIO} now implies that $W_t = e^{-ir^2/2t}e^{-itH}$ maps
scattering wavefront set to ordinary wavefront set, and vice versa.  In
particular, given the non-invariance of the Fourier transform, it is best
from a geometrical point of view to think of $W_t$ as a ``scattering Fourier integral operator''
with a wavefront relation interchanging $\RR^n \times S^{n-1}$ and $S^{n-1}
\times \RR^n.$

\section{The canonical relation}

The canonical relation of $W_t$ as a scattering FIO is related to the limit
of geodesic flow.  Given $(z,\hat \zeta) \in S^* \RR^n \cong \RR^n \times
S^{n-1}$ we let $\gamma(t)$ be the unit speed geodesic with $\gamma(0)=z$,
$(\gamma'(0))^i=g^{ij}\hat\zeta_j,$ and define (in the short range case\footnote{In the long range case there is a logarithmic divergence in $\lim t - |\gamma(t)|$ which needs to be removed. We omit the details.})
\begin{equation}
\begin{aligned}
S_f(z,\hat\zeta) &= (\theta, \xi) \in S^{n-1} \times \RR^n,\ \text{with } 
\xi=\lambda \theta+\mu,\\
\theta &= \lim_{t\to\infty} \frac{\gamma(t)}{|\gamma(t)|},\\
\lambda&=\lim_{t\to\infty} t-\abs{\gamma(t)}\\
\mu &= \lim_{t\to\infty} \abs{\gamma(t)} \left(\theta-\frac{\gamma(t)}{\abs{\gamma(t)}}\right).
\end{aligned}
\label{canonical}
\end{equation}
We also define $S_b(z, \hat \zeta) = -S_f(z, -\hat \zeta)$. These are the forward and backward sojourn relations, respectively. 
We interpret the components of these maps as follows: $\theta$ is the
asymptotic direction of the geodesic, $\lambda$ is the ``sojourn time,'' a
measure of how long the geodesic lingers in the finite part of $\RR^n$
before heading off to infinity, and $\mu$ measures the angle of contact of
the geodesic with the sphere at infinity, or, equivalently, distinguishes
different geodesics among the pencil of all geodesics with the same
asymptotic direction.

We may endow $\RR^n \times S^{n-1}$ with a contact structure using the
canonical one-form on $S^*\RR^n;$ by switching coordinates, as in the
definition of scattering wavefront set, we may thus endow $S^{n-1}\times
\RR^n$ with a contact structure.
\begin{proposition}
The maps
$$
S_f, S_b: \RR^n \times S^{n-1} \to S^{n-1} \times \RR^n
$$
are contact diffeomorphisms.
\end{proposition}
Hence $S_f$ and $S_b$ are eligible to be quantized to scattering Fourier integral
operators.

\begin{theorem}
The canonical relation of $W_t$ is $\abs{t}^{-1}S_f$ for $t<0$ and
$\abs{t}^{-1} S_b$ for $t > 0$, with the scaling acting in the fiber
variable.
\end{theorem}

\section{A sketch of the construction}

In \cite{Hassell-Wunsch1} a result considerably more detailed than
Theorem~\ref{par} is proven: it is shown that $W_t$ is a
\emph{scattering-fibered Legendrian distribution} in both space and time
variables simultaneously.  Legendrian distributions are a class of
distributions introduced by Melrose-Zworski in scattering theory
\cite{MR96k:58230}; they are given by oscillatory integrals and on $\RR^n$
are essentially the same as Fourier transforms of H\"ormander's Lagrangian
distributions.  The scattering-fibered Legendrians are a refinement
introduced by Hassell-Vasy \cite{MR2002i:58037,MR2001d:58034}.  This finer
characterization of $W_t$ amounts to the statement that it is a linear
combination of distributions of the form
\begin{equation}
t^{-\frac n2 -\frac k2} \int_{U\Subset \RR^k} a(t,r^{-1},\theta,z,v) e^{i
  \phi(r^{-1},\theta,z,v)r/t}\, dv
\label{phiform}\end{equation}
with $a$ and $\phi$ smooth in the short-range case (and having expansions
in the system $r^{-n} (\log r)^k$ for $k \leq n$ in general), and with
$\phi$ satisfying a non-degeneracy condition; see
\cite{MR2002i:58037,MR2001d:58034} for more details.  Note the special role
that $t$ plays in the phase.  All the geometric information is, as usual,
encoded in the phase function $\phi.$

The main step in the proof is the construction of a parametrix for the
propagator in the class of scattering-fibered Legendrians. One begins this
process for $t$ near zero near the diagonal of $\RR^n \times \RR^n$ (and
within the support of $\chi$). Here an ansatz for the propagator is the WKB
expression
\begin{equation}
t^{-n/2} e^{i\Phi(z,w)/t} \sum_{j \geq 0} t^j a_j(z, w).
\label{ansatz1}\end{equation}
Solving for the phase in the usual way gives the eikonal equation 
\begin{equation}
\Phi = \frac1{2} |\nabla_z \Phi|^2
\label{eikonal}\end{equation}
for $\Phi$, which has a smooth solution $\Phi = \h d_g(z,w)^2$, where $d_g$ is the distance with respect to the Riemannian metric $g$. This ansatz cannot be used for $z$ outside the injectivity radius of $w$, but a more complicated WKB-type ansatz 
\begin{equation}
t^{-n/2 - k/2} \int_{\RR^k} e^{i\Phi(z,w,v)/t} \sum_{j \geq 0} t^j a_j(z, w,v) \, dv
\label{ansatz2}\end{equation}
remains valid for $z$ in any given compact set. The main problem is to write down a suitable ansatz encapsulating the asymptotics as $|z| \to \infty$. 

In general, one solves a nonlinear first order PDE such as \eqref{eikonal} by the method of characteristics. In this case the characteristics are given by curves 
$$
s \mapsto (z=\gamma(s),\hat\zeta_i=g_{ij}(\gamma'(s))^j, t=0, \tau = s^2/2)
$$ where $\gamma$ is the arc-length parametrized geodesic starting from
$(w, \hat \xi)\in S_w^*\RR^n,$ and $\tau$ is the dual variable to $t$. The
set of all such geodesics sweeps out a submanifold $L$ of $\RR^n_z \times
\RR^n_\zeta \times \RR_\tau$ which remains smooth beyond the injectivity
radius. It is a Legendrian submanifold with respect to a naturally defined
contact structure on $\RR^n_z \times \RR^n_\zeta \times \RR_\tau$. If the
projection of $L$ to the $z$ variable is a diffeomorphism, then one can
write $\tau = \Phi(z)$ and this determines the phase function in
\eqref{ansatz1}. When this is not true, then one has to use the more
complicated expression \eqref{ansatz2}, where one requires $k$ extra variables
in order to write $\tau$ as a function of $z$ and $v$ if the null space of
the differential of the projection to $\RR^n_z$ has dimension $k$.

Along the geodesic, the form \eqref{metric} of the metric implies  (in the short-range case) that
\begin{equation}
r(\gamma(s)) = s + \Sigma(w, \hat\xi) + O(s^{-1}), \quad
\zeta(\gamma(s)) = \frac{\gamma(s)}{|\gamma(s)|} + O(s^{-1}).
\label{asymptoticflow}
\end{equation} Hence $\tau(\gamma(s)) = r^2/2 + r \Sigma(w, \hat\xi) + O(1)$. Recall
that we are interested in the phase of the propagator after an
$e^{ir^2/2t}$ factor has been stripped out.  Hence defining $\phi = \Phi -
r^2/2$, we find that $\phi \sim r\Sigma(x,\hat\xi)$ at infinity. Since
$\Phi$ (or more geometrically, the Legendre submanifold $L$) has good
asymptotics at infinity, we can hope to write down an ansatz for the
propagator which is accurate out to infinity and encodes the oscillations
at spatial infinity as well as as $t \to 0$. In fact, this turns out to be
possible in the class of fibered-scattering Legendre distributions and this
gives a parametrix which differs from the true propagator by a kernel which
is Schwartz in $z$ and $O(t^\infty)$, together with all its derivatives, as
$t \to 0$.

Moreover, in the same way that $\Phi$ in the ansatz \eqref{ansatz2} is
associated to a Legendre submanifold, so is $\phi(0, \cdot, \cdot, \cdot)$ in the expression \eqref{phiform} for
$W_t$. This is proved essentially by symplectic reduction. The Legendre submanifold associated to $\phi$ at $r= \infty$ then turns out to be the graph of the
transformation $S_f$, which is therefore a contact transformation.

The long-range case is not essentially different.

\

In the special case of flat $\RR^n,$ the sojourn time $\lambda$
along the geodesic emanating from $z$ in direction $\theta$ is
$$
\lim_{t\to\infty} t-\abs{z+t\theta}=-z\cdot \theta.
$$ Note that this is exactly the phase of $W_t$ in the Euclidean case, with
no oscillatory integral necessary.  This is typical of the following
special geometric case: if for each $w$ and $\theta$ ranging over a pair of
open sets, there exists a unique geodesic $\gamma$ with $\gamma(0)=w$,
  $\lim_{t\to\infty}
\gamma(t)/\abs{\gamma(t)}=\theta,$ then we can let $S(w,\theta)$ denote the
sojourn time along this geodesic (i.e.\ the function $\Sigma$ of
\eqref{asymptoticflow}, but parametrized by different variables).  Then
$S(w,\theta)$ parametrizes the Legendrian distribution, and we may write
the propagator for $t\geq 0$ and $w\in U \Subset \RR^n$ simply as
$$
t^{-n/2} a(t,w,r^{-1},\theta) e^{i S(w,\theta)/t}
$$ with $a$ a smooth function. The sojourn time $S(w,\theta)$ is closely related
to a sojourn time defined on pairs of points in $S^{n-1}$ by Guillemin
\cite{Gu} in studying the high-frequency asymptotics of the scattering
matrix.

\section{Propagation}

As a consequence of the characterization of $e^{-i \tilde r^2/2t} e^{-itH}$ as a
scattering FIO, where $\tilde r$ is given by \eqref{rt} ($\tilde r$ is just equal to $r$ in the short-range case), we can now state a propagation theorem, describing
precisely when and where singularities can appear in $\RR^n$:
\begin{theorem}
Let $\psi(t)$ be a solution to \eqref{scheqn}. Fix a time $t_0$. Let $(z, \hat \zeta) \in S^* \RR^n$, and
let
$$
(\theta, \xi) = S_b(z, \hat \zeta), \quad (\tilde \theta, \tilde \xi) = S_f(z, \hat \zeta).
$$
Then the following are equivalent:
\begin{align*}
(z, \hat \zeta) &\in \WF(\psi(t_0)); \\
(\theta, \frac{\xi}{t-t_0}) &\in \WFsc (e^{\, i\frac{\tilde r^2}{2(t-t_0)}} \psi(t)), \ t > t_0; \\
(\tilde \theta, \frac{\tilde \xi}{t_0 - t}) &\in \WFsc (e^{\, i\frac{\tilde r^2}{2(t-t_0)}} \psi(t)), \ t < t_0.
\end{align*}

\label{thm:propagation}
\end{theorem}

This is a refinement of certain propagation results of \cite{Wunsch1},
which were in terms of quadratic scattering wavefront set $\WFqsc$.  In particular,
one of the main results of \cite{Wunsch1} stated that if
$\theta=\pi(S_b(z,\hat\zeta))$ denotes the backwards limit of the geodesic
through $z,\hat\zeta$ then for $t>0,$
\begin{equation}
(\theta,-\theta/2t) \notin \WFqsc \psi_0 \Longrightarrow (z,\hat\zeta)
\notin \WF \psi(t).
\label{wfqsc}
\end{equation} This clearly follows from Theorem~\ref{thm:propagation} and the
following wavefront set computations:
\begin{proposition}
$$
(\theta,\xi) \in \WFqsc u \Longleftrightarrow (\theta,\xi-\alpha \theta/2)
  \in \WFqsc e^{-i \alpha r^2/2} u.
$$
\label{prop:wf1}
\end{proposition}
\begin{proposition}
If $(\theta,0) \notin \WFqsc u$ then $(\theta,\xi)\notin \WFsc u$ for any
finite $\xi.$
\label{prop:wf2}
\end{proposition}
The proofs of Propositions~\ref{prop:wf1}--\ref{prop:wf2}, which are not
difficult, are best carried out using the more sophisticated definitions of
wavefront set, involving associated calculi of pseudodifferential
operators.

\section{Examples}

We conclude with a pair of examples which exhibit some extremes of behavior
which solutions of \eqref{scheqn} may exhibit.

First, we consider the Schr\"odinger equation in flat $\RR^n$ with $V=0.$
Let $$\psi_0 (z) = (-2\pi i)^{-1/2} e^{-iz_1^2/2};$$ note that this is just
the one-dimensional fundamental solution, evaluated at time $t=-1,$ and
extended to be constant in the variables $z_2,\dots,z_n.$ The quadratic
scattering wavefront set of $\psi_0$ satisfies the hypotheses of
\eqref{wfqsc} over two points in $S^{n-1}_{\hat z},$ given by $(\pm 1, 0,
\dots, 0),$ hence \eqref{wfqsc} guarantees that the wavefront set of
$\psi(t)$ is confined at most to $\{t=1,\zeta=(\pm 1,0,\dots,0)\}.$ On the
other hand, we can compute exactly: $\psi(1,z) = \delta(z_1).$ Hence two
points in the qsc wavefront set are able to produce an entire hyperplane of
singularities.  By contrast, Theorem~\ref{thm:propagation} gives a more
satisfactory picture of the propagation phenomenon: we have
$$e^{ir^2/2}\psi_0=e^{i(z')^2/2}$$ (where $z'=(z_2,\dots,z_n)$).  Since the
Fourier transform of this function is just
$$
\frac{1}{i^{(n-1)/2} \sqrt{2\pi}} e^{-i(\zeta')^2/2} \delta(\zeta_1),
$$
we have
$$
\WFsc (e^{ir^2/2}\psi_0) = \{((\pm 1,0,\dots,0), (0,\zeta')): \zeta' \in\RR^{n-1}\}.
$$ This is pair of hyperplanes over two points in $S^{n-1}_{\hat \zeta}$
which maps \emph{diffeomorphically} to $\WF \psi(1,z)$ according to the sojourn
relation.

\

Another extreme case for the results discussed above occurs if we work on
$\RR^1$ with $V=0$, and take $$\psi_0(z) = e^{-iz^2/2} \Ai(z).$$ The
quadratic-scattering wavefront set of $\psi_0$ is just $(-1,1/2)$
(Intuitively speaking, the factor of $\Ai$ kills the wavefront set at $\hat
z = +1$ owing to its exponential decay in that direction, but has no effect
on the wavefront set of the $e^{-iz^2/2}$ factor at $\hat z=-1$ owing to
its slower oscillation.)  Hence \eqref{wfqsc} permits wavefront set only
at $t=1,$ in the direction $\zeta=+1.$ On the other hand, we may compute
exactly to find that
$$
\psi(z,1) = (-2\pi i)^{1/2} e^{i\left(\frac{z^3}3+\frac{z^2}2\right)} \in \mathcal{C}^\infty(\RR);
$$ hence no wavefront set appears after all.  This situation is accounted
for in the results of \cite{Wunsch1} by the use of a wavefront set that is
uniform in time, and in terms of which this solution is singular everywhere
on $\RR$ at $t=1.$ In terms of Theorem~\ref{thm:propagation}, which is
about the wavefront set of $\psi(\cdot, t)$ for each \emph{fixed} $t$, the
explanation is as follows: both $e^{ir^2/2} \psi_0$ and $\psi(z,1)$ have
scattering wavefront set only at the \emph{corner} $S^0 \times S^0,$
alluded to in the remark following Definition~\ref{defn:scwf}. As there
are no points $(\hat z, \zeta)$ in $\WFsc (e^{ir^2/2} \psi_0)$ with $\zeta$
finite, Theorem~\ref{thm:propagation} asserts that $\psi(1,z)$ is smooth.

\bibliography{all}
\bibliographystyle{amsplain}
\end{document}